\documentclass[12pt]{amsart}
\usepackage{amsmath, amsthm,amssymb, amscd, amsfonts}
\usepackage{fontsmpl}
\newtheorem{theorem}{Theorem}[section]
\newtheorem{corollary}[theorem]{Corollary}

\newtheorem{proposition}[theorem]{Proposition}
\newtheorem{lemma}[theorem]{Lemma}
\theoremstyle{definition}

\theoremstyle{remark}
\newtheorem{remark}[theorem]{Remark}

\numberwithin{equation}{section}

\newcommand\id{\operatorname{id}}

\newcommand\Tr{\operatorname{Tr }}

\newcommand\car{\operatorname{char}}
\newcommand\ord{\operatorname{ord}}
\newcommand\cop{\operatorname{cop}}
\newcommand\End{\operatorname{End}}
\newcommand\Alg{\operatorname{Alg}}

\newcommand\Rad{\operatorname{Rad}}
\newcommand\Map{\operatorname{Map}}
\newcommand\Perm{\operatorname{Perm}}
\newcommand\op{\operatorname{op}}
\newcommand{\eps}{\varepsilon}
\newcommand\Ker{\operatorname{Ker}}

\def\NN{\mathbb{N}}
\def\ZZ {\mathbb{Z}}
\def\QQ{\mathbb{Q}}

\def\GG{\mathbb{G}}

\begin{document}

%\rightline{}

\renewcommand{\theequation}{\arabic{equation}}

\title{On Hopf Algebras of dimension $p^{3}$}
\author[Gast\' on Andr\' es Garc\'\i a]{Gast\' on Andr\' es Garc\'\i a $^{*}$}
\address{Mathematisches Institut der Universit\" at M\" unchen,
\smallbreak
Theresienstr. 39, D-80333 M\" unchen, Deutschland.}
\email{garcia@mathematik.uni-muenchen.de}
\thanks{2000 {\it Mathemathics Subject Classification}: 16w30.
$^{*}$ Stipendiat des deutschen akademischen 
Austauschdienstes (DAAD), Germany.
This work will be part of the author's PhD thesis 
written under the supervision of Prof. N. Andruskiewitsch and
Prof. H.-J. Schneider.}

%\subjclass{}
%\date{\today.}

\begin{abstract}

We discuss some general results on 
finite-dimensional Hopf algebras over 
an algebraically closed field $k$ of characteristic zero
and then apply them to Hopf
algebras $H$ of dimension $p^{3}$ over $k$. There are 10 cases according to 
the group-like elements of $H$ and $H^{*}$. We show that 
in 8 of the 10 cases, it is 
possible to determine the structure of the Hopf algebra.
We also give a partial classification of the quasitriangular 
Hopf algebras of dimension $p^{3}$
over $k$, after studying extensions of a group algebra of 
order $p$ by a Taft algebra of dimension
$p^{2}$. In particular, we prove 
that every ribbon Hopf algebra of dimension $p^{3}$ 
over $k$ is either a group algebra
or a Frobenius-Lusztig kernel. Finally, using 
some results from \cite{AN} and \cite{BD} on bounds for the dimension 
of the first term $H_{1}$
in the coradical filtration of $H$, we give the complete classification 
of the quasitriangular Hopf algebras of dimension 27.
\end{abstract}

\maketitle

\section{Introduction}
We work over an algebraically closed field $k$ of
characteristic zero. Let $p$ be an odd prime number and let $\GG_{p}$ be
the cyclic group of $p$-th roots of unity.
We denote by 
$T(q)$, the Taft algebra of parameter $q \in \GG_{p} \smallsetminus \{1\}$, see Remark
\ref{taft} below. 
Hopf algebras of dimension 8 were classified by Williams \cite{Wi}; 
Masouka \cite{Mk2} and 
Stefan \cite{St} gave later a different proof of this fact.
In general, the classification problem of Hopf algebras over $k$ 
of dimension $p^{3}$ is still open. However,
the classification is known for semisimple or pointed
Hopf algebras.
Semisimple Hopf algebras of dimension $p^{3}$ were classified by Masuoka \cite{Mk1}; 
there
are $p+8$ isomorphism types, namely 

\begin{enumerate}
\item[(a)] Three group algebras of abelian groups.
\item[(b)] Two group algebras of non-abelian groups, and their duals.
\item[(c)] $p+1$ self-dual Hopf algebras which are neither commutative nor
cocommutative. They are extensions of 
$k[\ZZ/(p)\times \ZZ/(p)]$ by $k[\ZZ/(p)]$.
\end{enumerate}

Pointed non-semisimple Hopf algebras of dimension $p^3$ were classified
in \cite{AS2}, \cite{CD} and \cite{Stv}, by different methods.  
The explicit list is the following, where $q \in \GG_{p} \smallsetminus \{1\}$:
\begin{enumerate}
\item[(d)] The tensor-product Hopf algebra $T(q) \otimes k[\ZZ /(p)]$.
\item[(e)] $\widetilde{T(q)}  := k<g,\ x |\ gxg^{-1} =
q^{1/p}x,\ g^{p^{2}} = 1,\ x^{p} = 0>$ ($q^{1/p}$ a $p$-th
root of $q$), with comultiplication 
$\Delta (x) = x\otimes g^{p} + 1\otimes x,\ \Delta (g) = g\otimes g$.
\item[(f)] $\widehat{T(q)} := k<g,\ x |\ gxg^{-1} =
qx,\ g^{p^{2}} = 1,\ x^{p} = 0>$, with comultiplication
$\Delta (x) = x\otimes g + 1\otimes x,\ \Delta (g) = g\otimes g$.
\item[(g)] ${\bf r}(q) := k<g,\ x |\ gxg^{-1} = qx,\
g^{p^{2}} = 1,\ x^{p} = 1 - g^{p}>$, with comultiplication
$\Delta (x) = x\otimes g + 1\otimes x,\ \Delta (g) = g\otimes g$.

\smallbreak
\item[(h)] The Frobenius-Lusztig kernel ${\bf u}_{q}({\mathfrak {sl}}_{2}) 
:= k<g,\ x,\ y |\ gxg^{-1} = q^{2}x,$
$gyg^{-1} = q^{-2}y$, $g^{p} = 1$, $x^{p} = 0$,
$y^{p} = 0$, $xy - yx = g - g^{-1}>$, with comultiplication 
$\Delta (x) = x\otimes g + 1\otimes x,\ \Delta (y) =
y\otimes 1 + g^{-1}\otimes y,\
\Delta (g) = g\otimes g$.
\item[(i)] The book Hopf algebra ${\bf h}
(q, m) := k<g,\ x,\ y |\ gxg^{-1} = qx,\ gyg^{-1} = q^{m}y,\
g^{p} = 1,\ x^{p} = 0,\ y^{p} = 0,\ xy - yx = 0>,\ m \in \ZZ /(p) 
\smallsetminus \{0\}$,  
with comultiplication
$\Delta (x) = x\otimes g +
1\otimes x,\ \Delta (y) = y\otimes 1 + g^{m}\otimes y,\
\Delta (g) = g\otimes g$.
\end{enumerate}

Furthermore, there are two examples of non-semisimple but also non-pointed
Hopf algebras of dimension $p^{3}$, namely
\begin{enumerate}
\item[(j)] The dual of the Frobenius-Lusztig kernel, 
${\bf u}_{q}({\mathfrak {sl}}_{2})^{*}$.
\item[(k)] The dual of the case (g), ${\bf r}(q)^{*}$.
\end{enumerate}
 
There are no isomorphisms
between different Hopf algebras in the list. Moreover,
the Hopf algebras in cases (d), \ldots, (k) 
are not isomorphic for different values of 
$q \in \GG_{p} \smallsetminus \{1\}$, except for the book algebras,
where ${\bf h}
(q, m)$ is isomorphic to ${\bf h}
(q^{-m^{2}}, m^{-1})$.
In particular, the Hopf algebra $\widetilde{T(q)}$ does not depend, modulo
isomorphisms, upon the choice of the $p$-th root of $q$.  
The Hopf algebra ${\bf r}(q)$ was first
considered by Radford (see \cite{AS2}, \cite{R6}).

\smallbreak
We conjecture that any Hopf algebra $H$ of dimension $p^{3}$ is
semisimple or pointed or its dual is pointed, that is, $H$ is one
of the Hopf algebras of the list (a), \ldots, (k).

\smallbreak
In this paper we prove this conjecture under additional assumptions.
In Theorem \ref{dim1} and Corollary \ref{faltaba}
we show the simple modules of a crossed
product of a Taft algebra of dimension $p^{2}$ and a group algebra
of order $p$ are one-dimensional, that is the dual of the

crossed product is pointed.

\smallbreak
In Section \ref{generic} we discuss some general results on 
finite-dimensional Hopf algebras and then apply them to Hopf
algebras of dimension $p^{3}$. There are 10 cases according to 
the group-like elements of $H$ and $H^{*}$. We show that 
in 8 of the 10 cases, it is 
possible to determine the structure of the Hopf algebra.

\smallbreak
Let us say that a Hopf algebra $H$ of dimension $p^{3}$ is {\it strange} if $H$ is 
simple as a Hopf algebra, not semisimple and if $H$ and $H^{*}$ are not
pointed. It turns out that a Hopf algebra $H$ of dimension $p^{3}$ is 
isomorphic to a Hopf algebra of the list (a),\ \ldots,\ (k), or 
\begin{enumerate}
\item[(I)] $H$ is strange, $G(H) \simeq \ZZ / (p)\mbox{ and }G(H^{*}) = 1,$ or
\item[(II)] $H$ is strange, $G(H) \simeq \ZZ / (p),\mbox{ and }G(H^{*}) \simeq \ZZ / (p).$
\end{enumerate}
\noindent It is not known whether a strange Hopf algebra exists.

\smallbreak 
In the subsection \ref{pp}, we study non-semisimple 
Hopf algebras of dimension $p^{3}$ with $G(H) \simeq G(H^{*}) \simeq \ZZ /(p)$.
The order of the antipode of such a Hopf algebra is necessarily $2p$ or $4p$. If the
order is $2p$, then $H$ is a bosonization of the group algebra $k[G(H)]$. In this case
we believe that $H$ is isomorphic to a book Hopf algebra ${\bf h}(q, m)$, for some
$q\in \GG_{p} \smallsetminus \{1\}$, and $m\in \ZZ /(p) \smallsetminus \{0\}$.
If the order is $4p$, then $H$ satisfies (II), and all skew primitive elements of $H$ 
are trivial, that is, contained in $k[G(H)]$.

\smallbreak 

Radford and Schneider \cite{RS2}
conjectured that the square of the antipode of any 
finite-dimensional Hopf algebra
must satisfy a certain condition, which they called the
strong vanishing trace condition. If 
$H$ is a finite-dimensional Hopf algebra and $B$ is the
unique maximal semisimple Hopf subalgebra of $H$, then it follows from the
conjecture that
the order of the square of the antipode of $H$ must divide
$\dim H / \dim B$, see \cite[Thm. 6]{RS2}. 
In particular, if the dimension of $H$ is $p^{3}$ and $|G(H)| = |G(H^{*})| = p$ or 
$|G(H)| = p$ and $|G(H^{*})| = 1$, then the order of the antipode should be $2p$.

\smallbreak
It is well-known that  
the Frobenius-Lusztig kernels ${\bf u}_{q}({\mathfrak {sl}}_{2})$ and the 
group algebras admit a quasitriangular structure (see
e. g. \cite[IX. 7]{K}). We prove in Section \ref{qt} that
these two are the only quasitriangular Hopf algebras from the list
above. 
We also prove in Theorem \ref{principal} that there is no 
quasitriangular Hopf algebra of dimension $p^{3}$ 
which satisfies condition (I).
Namely, if $H$ is a quasitriangular Hopf algebra of dimension $p^{3}$, then
\begin{enumerate}
\item[$(i)$] $H$ is isomorphic to a 
group algebra or to ${\bf u}_{q}({\mathfrak {sl}}_{2})$, or
\item [$(ii)$] $H$ satisfies (II) and the map $f_{R}:H^{* cop} \to H$ 
is an isomorphism. Moreover,
$H$ and $H^{*}$ are minimal quasitriangular, $1 = <\beta, x>$, for all
$\beta \in G(H^{*})$, $x\in G(H)$ and $\ord {\mathcal S} = 4p$.
\end{enumerate}

\smallbreak 
As a consequence, we show in Corollary \ref{ribbon}
that every ribbon Hopf algebra of dimension
$p^{3}$ is either a group algebra or a Frobenius-Lusztig kernel.
\smallbreak
Finally, using 
some results from \cite{AN} and \cite{BD} on bounds for the dimension 
of the first term $H_{1}$
in the coradical filtration of $H$, we classify 
quasitriangular Hopf algebras of dimension 27.

\subsection{Conventions and Preliminaries}\label{conv}
Our references for the theory of Hopf algebras are \cite{Mo}, \cite{K}, 
\cite{Sch} and \cite{Sw}. The antipode of a Hopf algebra $H$ is
denoted by $\mathcal  S$, or ${\mathcal  S}_H$ if special emphasis
is needed. The Sweedler notation is used for the
comultiplication of $H$ but dropping the summation symbol. The
group of group-like elements of a coalgebra $C$ is denoted by $G(C)$.
The modular group-like elements $g\in H$ and $\alpha \in H^{*}$ are
defined by 
$$ \Lambda x = <\alpha,x>\Lambda ,\mbox{ for all }x\in H,\mbox{ and }
\beta \lambda = <\beta, g> \lambda,\mbox{ for all }\beta\in H^{*},$$
\noindent where $\Lambda \in H$ is a non-zero left integral and
$\lambda \in H^{*}$ is a non-zero right integral. 
There is a 
formula for ${\mathcal S}^{4}$ in terms of $\alpha$ and $g$:  
$${\mathcal S}^{4}(h) = g(\alpha \rightharpoonup h \leftharpoonup 
\alpha^{-1}) g^{-1},\mbox{ for all }h \in H,$$
where $H^{*}$ acts on $H$ on the left by
$\beta \rightharpoonup h = h_{(1)}\beta (h_{(2)})$ and on the right by
$h\leftharpoonup \beta  = \beta(h_{(1)})h_{(2)}$, for all 
$\beta \in H^{*}$ and $h\in H$. Moreover, if $\lambda$ and $\Lambda$ are
such that $<\lambda, \Lambda>= 1$, then there are
formulas for the trace of any linear endomorphism $f$ on $H$: 
$$\Tr f = <\lambda, {\mathcal S}(\Lambda_{(1)})f(\Lambda_{(2)})>
= <\lambda, ({\mathcal S}\circ f)(\Lambda_{(1)})\Lambda_{(2)}>.$$
The formulas above are due to many authors, including Radford 
(see \cite{Mo}, \cite{R1}, \cite{Sch}).
Let $C$ be a coalgebra and $a,\ b \in G(C)$. The set of 
$(a, b)$-skew primitive elements of $C$ is defined by 
$$P_{a,b} = \{ c \in C |\ \Delta (c) = a\otimes c + c\otimes b\};$$

\noindent in particular, $k(a-b) \subseteq P_{a,b}$. We say that a skew primitive element
$c \in C$ is trivial if $c \in k[G(C)]$.

\smallbreak
A Hopf algebra $H$ is called {\it simple} if it does not contain any proper
normal Hopf subalgebra in the sense of \cite[3.4.1]{Mo}; $H$ is
called {\it semisimple}, respectively {\it cosemisimple}, if it
is semisimple as an algebra, respectively if it is cosemisimple
as a coalgebra. The sum of all simple subcoalgebras is called the
{\it{coradical}} of $H$ and it is denoted by $H_{0}$. If all
simple subcoalgebras of $H$ are one-dimensional, then $H$ is
called {\it pointed} and $H_{0} = k[G(H)]$. 
A finite-dimensional Hopf algebra $H$ is pointed if and only if 
all the simple $H^{*}$-modules are one-dimensional.
We also consider
the coradical filtration $H_{0} \subset H_{1} \subset \cdots$, of $H$ 
(see \cite[Chapter 5]{Mo}, \cite[Chapter IX]{Sw}).

\smallbreak

Let $H$ be a finite-dimensional Hopf algebra over $k$, then the following
are equivalent (see \cite{LR}, \cite{LR2}, \cite[Prop. 2]{R3} and
\cite[Cor. 3.5]{Sch}):
\begin{enumerate}
\item[$(a)$] $H$ is semisimple, 
\item[$(b)$] $H$ is cosemisimple, 
\item[$(c)$] ${\mathcal S}^2 = \id$, 
\item[$(d)$] $\Tr {\mathcal S}^2 \neq 0$,
\end{enumerate}

\noindent where $\Tr$ denote the trace map.

\begin{remark}\cite{Z} \label{oddim}
If $H$ is an odd-dimensional Hopf algebra and ${\mathcal  S}^4 =
\id$, then $H$ is semisimple. Therefore, if $H$ is a
non-semisimple Hopf algebra of odd dimension, either $G(H)$ or
$G(H^{*})$ is non-trivial.
\end{remark}

\begin{remark}\label{taft}
Let $N \geq 2$ be an integer and let  $q \in k$  be
a primitive $N$-th root of unity.  Recall that the {\it Taft
algebra} $T(q)$ is the $k$-algebra presented by generators $g$ and
$x$ with relations $g^N = 1$, $x^N = 0$ and $gx = q xg$; $T(q)$
carries a Hopf algebra structure, determined by
$$ \Delta g = g \otimes g, \qquad \Delta x = x \otimes 1 + g
\otimes x.$$
\noindent Then $
\eps (g) = 1,\ \eps(x) = 0,\ 
\mathcal S(g) = g^{-1}$, and $\mathcal S(x) = -g^{-1}x$.
It is known that $T(q)$ is a pointed Hopf algebra, with $G(T(q)) =
\langle g  \rangle \simeq {\mathbb Z} / (N)$. The proper Hopf
subalgebras of $T(q)$ are contained in $k \langle g \rangle$,
whence they are semisimple. We also have: 
\begin{enumerate}
\item[$(i)$] $T(q) \simeq T(q)^*$,
\item[$(ii)$] $T(q) \simeq T(q')$ if and only if $q = q'$. 
\end{enumerate}
\noindent It is not difficult
to see that ${T(q)^*}^{\cop} \simeq T(q)^{\op} \simeq T(q^{-1})$.

The square of the antipode of $T(q)$ coincides with the inner
automorphism induced by $g$. Therefore, ${\mathcal S}^4 \neq \id$
if $N > 2$.
\end{remark}

\section{Extensions of a Taft algebra by the group algebra of order $p$}
\smallbreak 
In the list of Hopf algebras of dimension
$p^3$ given above, the cases (d), (f) and (g) are extensions of
Taft algebras by the group algebra of order $p$; that is, they fit into an 
exact sequence of finite-dimensional Hopf algebras 
$$
1\to k[\ZZ/(p)] \xrightarrow{\iota} H \xrightarrow{\pi}
T(q) \to 1.
$$
\smallbreak
We show that in some sense the converse is also true, that is, 
if $H$ is an extension of a Taft algebra
by the group algebra of order $p$, then $H$ is necessarily pointed. 
\smallbreak
We recall some definitions and formulas 
for cleft extensions. For examples and a
characterization of these extensions see
\cite[Chapter 7]{Mo}. 

By \cite[Thm. 2.4]{Sch2}, every extension 
$1\to A \xrightarrow{\iota} H \xrightarrow{\pi}
B \to 1,$
of finite-dimensional Hopf algebras is cleft; {\it{i.e.}} there exists 
a right $B$-comodule map $\gamma: B \to H$, which is convolution invertible 
and preserve the unit and the counit.  
Using this map, one can construct a weak action
of $B$ on $A$ and a convolution invertible 2-cocycle 
$\sigma \in$ Hom $(B\otimes B,A)$,
that give $H$ the structure of a crossed product $A\#_{\sigma} B$ 
of $A$ with $B$. As vector spaces $A\#_{\sigma} B = A\otimes B$.
The weak action of $B$ on $A$ is defined by
$$
b\cdot a = \gamma (b_{(1)}) a \gamma ^{-1}(b_{(2)}),$$
\noindent for all $a \in A,\ b \in B$, and satisfies
$$
b\cdot (aa') = (b_{(1)}\cdot a)(b_{(2)}\cdot a') \mbox{ and }\ b\cdot 1 =
\eps (b) 1_{A},
$$
\noindent for all $b \in B,\ a,\ a' \in A$. 
The convolution invertible 2-cocycle $\sigma$ is
given by 
$$
\sigma (b,c) = \gamma (b_{(1)})\gamma (c_{(1)})\gamma^{-1} (b_{(2)}c_{(2)}),
$$
\noindent for all $b,\ c \in B$, and satisfies
$$
\begin{array}{c}
[b_{(1)}\cdot \sigma(b'_{(1)},\ b''_{(1)})]\sigma(b_{(2)},\ b'_{(2)}b''_{(2)}) 

=
\sigma(b_{(1)},\ b'_{(1)})\sigma(b_{(2)}b'_{(2)},\ b'')\\
\\
\mbox{ and}\qquad
\sigma(b,\ 1)= \sigma(1,\ b)= \eps (b)1,
\end{array}
$$
\noindent for all $b,\ b',\ b'' \in B$. The multiplication
on the crossed product $A\#_{\sigma} B$ is given by 
$$
(a\# b)(c\# d) = a (b_{(1)}\cdot c)\sigma (b_{(2)},d_{(1)})\# b_{(3)}d_{(2)},
$$
\noindent for all $a,\ c \in A,\ b,\ d \in B$. The unit in $A\#_{\sigma} H$
is $1\#1$. Here we write $a\#h$ for $a\otimes h$ as an element in
$A\#_{\sigma} B$.  

\smallbreak
If $B = k[\Gamma]$ is a group algebra of a group $\Gamma$, then 
$A *_{\sigma} \Gamma := A\#_{\sigma} k[\Gamma]$ is called the
$\Gamma$-crossed product of $\Gamma$ over $A$.
Observe that $A*_{\sigma} \Gamma$ is a $\Gamma$-graded algebra. 
Moreover, it is easy to see that $A*_{\sigma} \Gamma$ is 
characterized as the $\Gamma$-graded algebra which contains 
an invertible element in each component and whose 1-component
is $A$. In this case, we say that $A$ is the neutral component
of $A*_{\sigma} \Gamma$.

\begin{theorem}\label{dim1}
Let $H$ be a finite-dimensional Hopf algebra which fits into 
an extension of the form  
\begin{equation}\label{teo2.1} 
1\to A \xrightarrow{\iota} H \xrightarrow{\pi}
k[\Gamma] \to 1,
\end{equation}

\noindent where $\Gamma = \ZZ / (p)$ is the group algebra
of order $p$, $A^{*}$ is pointed and
the group of group-likes of $A^{*}$ has order $|G(A^{*}| \leq p$. 
Then $H^{*}$ is pointed.
\end{theorem}

\begin{proof}
To prove that $H^{*}$ is pointed, we show that every 
simple $H$-module is one-dimensional. To see this  
it is enough to show that the algebra $H / \Rad H$ is
commutative. Indeed, we have that 
$H / \Rad H$ is a semisimple algebra. If
it is commutative, then by the Artin-Wedderburn theorem it 
follows that every simple $H / \Rad H$-module, and
therefore every simple $H$-module, is one-dimensional. 

\smallbreak
Hence, we prove actually the following. Let 
$H = A*_{\sigma} \Gamma$ be a finite-dimensional 
$\Gamma$-crossed product with neutral component $A$
and suppose that
\begin{enumerate}
\item[(i)] $A /\Rad A \simeq \Map (X,k)$, and 
\item[(ii)] there exists an epimorphism $\pi : H \to k[\Gamma]$ of
$\Gamma$-graded algebras,
\end{enumerate}

\noindent where $\Map (X, k)$ is the set of functions 
on the finite set $X$ and $|X|\leq p$. 
Then $H / \Rad H$ is commutative. 

\smallbreak
Let $\overline{A} = A /\Rad A$.
$X$ can be naturally identified with the set 
$\{\delta_{0}, \ldots , \delta_{m} \}$ of 
primitive idempotents of $\overline{A}$. 

\smallbreak
Since for every $g \in \Gamma$, the map $r(g): A \to A$ given
by $r(g)(a) = g\cdot a$ for all $a\in A$ is an algebra map,
the radical $\Rad A$ of $A$ is stable under the weak action
of $\Gamma$ on $A$ and whence $(\Rad A) * \Gamma$ is a 
$\Gamma$-graded ideal in $H$, which is nilpotent. Therefore
we have a quotient algebra $\overline{H} = \overline{A} * \Gamma$,
which is also a $\Gamma$-crossed product. Since $\car\ k= 0$, 
$\overline{H}$ is semisimple, whence
$\overline{H} = H /\Rad H$. Moreover,
since $k[\Gamma]$ is a semisimple algebra, it follows that
$\pi (\Rad H) = 0$ and $\pi$ factorizes through $\overline{H}$. 
Denote by 
$\bar{\pi} : \overline{H} \to k[\Gamma]$ the map induced by this 
factorization.

\smallbreak
Let $\delta$ be a primitive idempotent of $\overline{A}$. Then
for all $g \in \Gamma$,
$g\cdot \delta$ is also a primitive idempotent of $\overline{A}$. Hence,
the weak action of $\Gamma$ associated to $\overline{H}$ arises from a 
group homomorphism, say $\alpha : \Gamma \to \Perm (X)$. 
Since for all $g \in \Gamma$ and $\delta_{i} \in X$, 
$$\bar{\pi} (\alpha(g)(\delta_{i}) *1) = \bar{\pi} (g\cdot \delta_{i} *1)=
\bar{\pi} ((1*g)(\delta_{i} *1)(1*g^{-1}))= 
\eps (\delta_{i})1 = \delta_{i,0},$$

\noindent it follows that 
$\Gamma$ must fix the unique element $\delta_{0}$ in $X$ that
does not vanish under $\bar{\pi}$. Since $\Gamma = \ZZ /(p)$ 
and $p$ does not divide $(|X| -1)!$, the homomorphism $\alpha$
must be trivial. This implies that the weak action of $\Gamma$ on
$\overline{A}$ is trivial and therefore that the $\Gamma$-crossed product 
$\overline{H}$ is commutative.
\end{proof}

\begin{corollary}\label{faltaba}
Let $H$ be a Hopf algebra over $k$ of dimension $p^{3}$ which is an extension
of a Taft algebra $T(q)$ by a group algebra $k[\ZZ/(p)]$, 
that is, $H$ fits into the
following 
exact sequence of finite-dimensional 
Hopf algebras: 
\begin{equation}\label{T1}
1\to k[\ZZ/(p)] \xrightarrow{\iota} H \xrightarrow{\pi}
T(q) \to 1.
\end{equation}
Then $H$ is pointed.
\end{corollary}

\begin{proof} If we dualize the sequence (\ref{T1}), 
we see that $H^{*}$ is an extension of a group algebra of order $p$
by a Taft algebra of order $p^{2}$,
$$
1\to T(q) \xrightarrow{\iota^{*}} H^{*} \xrightarrow{\pi^{*}}
k[\ZZ/(p)] \to 1.
$$

Since the Taft algebra satisfies the hypothesis on $A$ in Theorem \ref{dim1},
with $n =1$ and $N_{1} =p$,  
by Theorem \ref{dim1} the simple $H^{*}$-modules are one-dimensional, 
which implies that $H$ is a pointed Hopf algebra.
\end{proof}

The following corollary states that if there exist a Hopf 
algebra $H$ of dimension $p^3$ which is not isomorphic to a Hopf
algebra of the list (a),\ \ldots,\ (k), then $H$
is necessarily strange.

\begin{corollary}\label{simplicidad}
Let $H$ be a non-semisimple Hopf algebra of dimension $p^{3}$ such that 
$H$ and $H^{*}$ are both non-pointed. Then $H$
is simple as a Hopf algebra.
\end{corollary}

\begin{proof}
Suppose that $H$ is not simple, then it contains a proper 
Hopf subalgebra $K$ which is normal and non-trivial.
Hence, we have an extension of Hopf algebras
\begin{equation}\label{ext1}
1\to K \xrightarrow{\iota} H \xrightarrow{\pi} 
D \to 1, 
\end{equation}

where $D= H/K^{+}H$. Then it follows that $p^{3} = \dim K \dim D$, by Nichols-Zoeller,  
and the dimension of $K$ is $p$ or $p^{2}$.

\smallbreak
If dim $K = p^{2}$, then dim $D =p$. Moreover, by \cite[Thm. 2]{Z},
$D \simeq k[\ZZ/(p)]$ and by \cite[Thm. 5.5]{Ng}, 
$K \simeq T(q)$, a Taft algebra,
since $H$ is non-semisimple by assumption. Hence $H$ 
is an extension of a group algebra
by a Taft algebra and by the Corollary \ref{faltaba} , 
$H^{*}$ must be pointed, which is a contradiction.

\smallbreak
If dim $K= p$, applying the same argument to the dual extension of 
(\ref{ext1}) we get that $H$ is pointed,
which is also a contradiction.
\end{proof}

\bigbreak
\section{On Hopf algebras of dimension $p^3$}\label{generic}

In this section we first discuss some general results on
Hopf algebras of finite dimension and then we apply them to Hopf algebras
of dimension $p^3$. 

\subsection{General results}
The order of the antipode plays a very important r$\hat{\mbox{o}}$le
in the theory of finite-dimensional Hopf algebras. 
The following linear algebra lemma will help us to determine 
this order in some particular cases. Note that part $(c)$  
generalizes $(b)$.

\begin{lemma}\label{gral}
(a) \cite[Lemma 2.6, (i)]{AS1} 
Let $T$ be a linear transformation of
a finite-dimensional vector space $V \neq 0$ such that $\Tr T = 0$ and $T^{p} =\id$.
Let $q \in \GG_{p}\smallsetminus \{1\}$ and let $V(i)$ be the eigenspace
of $T$ of eigenvalue $q^{i}$, then $\dim V(i)$ is constant; in particular
$p$ divides $\dim V$.
\smallbreak
(b) \cite[Lemma 2.6, (ii)]{AS1} If $V$ is a vector space of dimension $p$
and $ T \in \End V$ satisfies 
$T^{2p} = \id,$ and $\Tr T = 0,$ then 
$T^{p} = \pm \id$.
\smallbreak 
(c) Let $n \in \NN$ and let $\omega$ be a $p^{n}$-th primitive 
root of unity in $k$.  
If $V$ is a vector space of dimension $p$
and $T \in \End V$ satisfies 
$T^{2p^{n}} = \id,$ and $\Tr T = 0,$ then there 
exists $m,\ 0\leq m\leq p^{n-1}-1$ such that $T^{p} = \pm \omega^{mp}\id$.
\end{lemma}

\begin{proof}
To prove $(c)$ we follow the approach of \cite{AS1}. The crucial point here
is that the minimal polynomial over $\QQ$ of a $p^{n}$-th root of unity in $k$ is 
known and that $V$ is $p$-dimensional. 

\smallbreak
Since $T^{2p^{n}} = \id,$ the eigenvalues of $T$ are 
of the form $(-1)^{a}\omega ^{i},$ where $a \in \{0,1\}$ and $0\leq i\leq p^{n}-1$.   
Let $V_{a,i} := \{ v \in V : T (v) = (-1)^{a}\omega ^{i}v\}$ be the eigenspace
of $T$ of eigenvalue $(-1)^{a}\omega ^{i}.$ Since $\Tr T = 0$, we have 
\begin{equation}\label{polin}
 0 = \Tr T = \sum_{i=0}^{p^{n}-1} 
(\dim V_{0,i} - \dim V_{1,i})\omega ^{i}.
\end{equation}

\noindent Define now $P \in \ZZ[X]$ by 
$P(X) = \sum_{i=0}^{p^{n}-1} (\dim V_{0,i} - \dim V_{1,i})X^{i}.$
Then deg $P \leq p^{n} -1,\ P(\omega) = 0$ and the 
number of coefficients of $P$ different from
zero is less or equal than $p$, since $V = \bigoplus_{a,i} V_{a,i}$ and dim $V = p$. 
Moreover, the
minimal polynomial $\Phi_{\omega}$ of $\omega$ 
over $\QQ$ must divide $P$, since $\omega$ is a root of $P$. Hence there exists
$Q \in \ZZ[X]$ such that $P = Q\Phi_{\omega}$, where $Q$ is the zero
polynomial or  
deg $Q \leq p^{n-1}-1$, since deg $\Phi_{\omega}
= \varphi (p^{n}) = p^{n-1}(p-1)$ 
(where $\varphi$ is the Euler's function). 

\smallbreak
Define $V_{+} = \bigoplus_{i= 0}^{p^{n}-1} V_{0,i}$ and 
$V_{-} = \bigoplus_{i= 0}^{p^{n}-1} V_{1,i}$, then it is clear that 
$V = V_{+} \oplus V_{-}$.
If $Q$ is the zero polynomial, then $P$ is also the zero polynomial and it  
follows that 
$\dim V_{0,i} = \dim V_{1,i},$ for all $0\leq i\leq p^{n}-1$.
But this implies that dim $V_{+}$ = dim $V_{-}$, from which follows that
the
dimension of $V$ is even, a contradiction.

\smallbreak
Therefore we can assume that $Q$ is not the zero polynomial.
Suppose that $Q(X) = \sum_{j=0}^{p^{n-1}-1} d_{j}X^{j}$, where
$d_{m}\neq 0$ for
some $0\leq m\leq p^{n-1}-1$, and  
recall that 
$\Phi_{\omega}(X) = X^{p^{n-1}(p-1)} + X^{p^{n-1}(p-2)} + \cdots + X^{p^{n-1}} + 1.$
Then we have that
$$
P(X)  = \left( \sum_{i=0}^{p^{n-1}-1} d_{i}X^{i}\right)
\left(\sum_{j=0}^{p-1} (X^{p^{n-1}})^{j} \right) 
= \sum_{i=0}^{p^{n-1}-1} d_{i}\left(\sum_{j=0}^{p-1} X^{i + jp^{n-1}}\right).$$

Since the number of non-zero coefficients of $P$ is not zero 
and less or equal than $p$, there exists a unique 
$l,\ 0\leq l\leq p^{n-1}-1$ such that $d_{l}\neq 0$ and
$d_{i} = 0$ for all $i\neq l,\ 0\leq i\leq p^{n-1}-1$.
Hence, $l = m$ and then
$P(X) = d_{m}X^{m}\Phi_{\omega}$, implying for all 
$0\leq j\leq p-1$, $0\leq i\leq p^{n}-1$ such that
$i\neq m + jp^{n-1}$ that
$$
%\begin{array}{c} 
\dim V_{0,m+ jp^{n-1}} - \dim V_{1,m+ jp^{n-1}}  = d_{m},
\mbox{ and }\dim V_{0,i} - \dim V_{1,i} = 0.  
%\end{array}
$$

\noindent From these equalities it follows that 
dim $V_{+} -$ dim $V_{-} = d_{m}p$; 
as dim $V_{+}$ + dim $V_{-} = p$,
one concludes that 
$V_{+} = \bigoplus_{j=0}^{p-1} V_{0, m+ jp^{n-1}} = V$
or $V_{-} = \bigoplus_{j=0}^{p-1} V_{1, m+ jp^{n-1}} = V$, and this implies
that $T^{p}(v) = \pm (\omega^{m})^{p}v= \pm \omega^{m p} v$, for all $v \in V$.
\end{proof}

\begin{proposition} \cite[Prop. 5.1]{AS1} \label{esta}
Let $H$ be a finite-dimensional Hopf algebra whose antipode
${\mathcal  S}$ has order $2p$. Assume also that $H$ contains a
cosemisimple Hopf subalgebra $B$ such that $\dim H= p\
\dim B$. Then $B$ is the coradical of $H$.
\qed \end{proposition}

We use Lemma \ref{gral} $(c)$ to improve a result due to 
Radford and Schneider \cite{RS}. As shown in the remark below,
this result gives an alternative proof of the
classification of Hopf algebras of dimension $p^{2}$. 

\begin{proposition} \label{schrad}
Let $H$ be a finite-dimensional non-semisimple Hopf algebra which 
contains a commutative Hopf subalgebra $B$ such that $\dim
H = p\ \dim B$. If ${\mathcal  S}^{4p^{n}} = \id,$ for some $n \in \NN$,
then ${\mathcal S}^{2p} = \id$, and consequently 
$B$ is the coradical of $H$.
\end{proposition}

\begin{proof}
We follow the proof of \cite{RS}.
By assumption on $k$, we know that $B$ is semisimple, since it is
commutative. Let $\{ e_{j}\}_{1\leq j \leq s}$, be the central
primitive idempotents of $B$ such that
$ 1_{B} = e_{1} + \cdots + e_{s}.$
Let $I_{j} = H e_{j}$; then we can write $H = \bigoplus_{j=1}^{s}
I_{j}$.

\smallbreak {\noindent{\bf Claim 1.}}  $\Tr ({\mathcal
S}^{2}|_{I_{j}}) = 0,\ \forall\ 1\leq j \leq s$.

\smallbreak
It is clear that ${\mathcal  S}^{2}(I_{j}) \subseteq I_{j}$, since
$${\mathcal  S}^{2}(I_{j})  = {\mathcal  S}^{2}(H e_{j})  =
{\mathcal  S}^{2}(H){\mathcal  S}^{2}(e_{j}) =  
{\mathcal  S}^{2}(H)(e_{j}) \subseteq H e_{j} = I_{j}.$$

Let $P_{j} : H \to I_{j}$ be the linear projection given by
$P_{j}(a) = a e_{j},$ for all $a\ \in H$, and $r(h)$ the
linear map given by $r(h)(x) = xh, \forall\ x,\ h \in H$. Then
${\mathcal S}^{2}|_{I_{j}} = {\mathcal  S}^{2}\circ P_{j} =
{\mathcal S}^{2}\circ r(e_{j})$, since these maps coincide
on $I_{j}$. Let $\lambda \in H^{*}$ be a right integral and
$\Lambda \in H$ a left integral such that $<\lambda ,\ \Lambda
> = 1$; then $(\Lambda_{(1)} ,\ {\mathcal  S}(\Lambda_{(2)}))$ are
dual bases for the Frobenius homomorphism $\lambda$, (see
\cite{Sch}). By Radford's trace formula,
we get that
$$\begin{array}{ll}
\Tr ({\mathcal  S}^{2}|_{I_{j}}) & =  \Tr (({\mathcal  S}^{2})\circ r(e_{j}))
=  <\lambda ,\ {\mathcal  S}(\Lambda_{(2)})
[{\mathcal  S}^{2}\circ r(e_{j})] (\Lambda_{(1)}) > \\
\\
& =  <\lambda ,\ {\mathcal  S}(\Lambda_{(2)})\ 
{\mathcal  S}^{2} (\Lambda_{(1)} e_{j}) >\\
\\
& =  <\lambda ,\ {\mathcal  S}(\Lambda_{(2)}){\mathcal  S}^{2}
( \Lambda_{(1)}){\mathcal  S}_{B}^{2}( e_{j}) > \\
\\
& =  <\lambda ,\ {\mathcal  S}({\mathcal  S}(\Lambda_{(1)}) 
\Lambda_{(2)})\ e_{j} >
 =  <\lambda ,\ e_{j} ><\eps ,\ \Lambda > =  0,
\end{array}$$

\noindent{where} the last equality follows from the fact that 
$H$ is non-semisimple.

\smallbreak {\noindent{\bf Claim 2.}} $\dim I_{j} = p,\
1\leq j \leq s$.

\smallbreak
Since for all idempotents $e \in B,\ e\neq 0$, $He \cong H\otimes_{B} Be$, 
as vector spaces over $k$, we have that $\dim He = 
\dim_{B}\ H\  \dim Be$. But this dimension is $p$,
since $\dim Be  = 1$, because $B$ is commutative.

\smallbreak Let $T_{j} = {\mathcal  S}^{2}|_{I_{j}},\ 1\leq j \leq
s$. By the claims above, it follows that $\Tr (T_{j}) = 0,\
T_{j}^{2p^{n}} = \id_{I_{j}}$ and dim $I_{j} = p$, which implies by 
Lemma \ref{gral} $(c)$ that $T_{j}^{p} = \pm \omega_{j}^{m_{j}p}\id_{I_{j}}$, 
where $\omega_{j}$ a $p^{n}$-th root of unity and 
$0\leq m_{j} \leq p^{n-1}-1.$ 
Since $T_{j}(e_{j})
= {\mathcal  S}^{2}|_{I_{j}}(e_{j}) = e_{j},$ 
it follows that $m_{j} = 0$ 
and $T_{j}^{p} = \id_{I_{j}},\ \forall\ 1\leq j
\leq s$, that is,
${\mathcal  S}^{2p}|_{I_{j}} = \id_{I_{j}}$, for all
$1\leq j\leq s$. Finally, the proposition follows  
in view of Proposition \ref{esta}.
\end{proof}

\begin{remark}\cite{RS}
Radford and Schneider proved Proposition \ref{schrad} in the case 
where $n=1$ using Lemma \ref{gral} $(b)$. This result provides 
an alternative proof of the classification of Hopf algebras of
dimension $p^{2}$, which was recently finished by Ng \cite{Ng}.
Indeed, if $H$ is
semisimple, then by a result of Masuoka \cite{Mk3}, it is
isomorphic to a group algebra of order $p^2$. Now suppose that $H$
is non-semisimple. By Remark \ref{oddim}, $G(H)$ or $G(H)^{*}$ is
not trivial and the order of each group is less than $p^2$. 
Assume that $G(H)$ is not trivial. Then $\dim H$ $= p
|G(H)|$ 
and it follows from Radford's formula for ${\mathcal S}^{4}$ that
${\mathcal S}^{4p} = \id$. Hence, by the
proposition above, $H$ must be pointed and by \cite[Thm. A
(ii)]{AS1} $H$ is isomorphic to a Taft algebra.
\end{remark}

\bigbreak
Let $H$ be a finite-dimensional Hopf algebra and suppose that $G(H)$ is
abelian and its order is $p^{n},\ n\geq 0$. Following the
classification of finite abelian groups, we say that $G(H)$ is of type
$(p^{i_{1}},p^{i_{2}}, \ldots , p^{i_{s}})$ if 
$G(H) \simeq \ZZ /(p^{i_{1}}) \times \ldots \times \ZZ/ (p^{i_{s}})$.
Suppose in addition that $G(H^{*})$ is also abelian
and its order is a power of $p$. We say 
that the Hopf algebra $H$ is of {\it type}
$(p^{i_{1}},\ldots , p^{i_{s}} ; p^{j_{1}},\ldots , p^{j_{t}})$
if $G(H)$ and $G(H^{*})$ are of type $(p^{i_{1}},\ldots , p^{i_{s}})$ 
and $(p^{j_{1}},\ldots , p^{j_{t}})$, respectively.  

\smallbreak
If $H$ is a non-semisimple
Hopf algebra of dimension $p^{3}$, then $G(H)$ and $G(H^{*})$ are abelian 
and their orders are powers of $p$, by Nichols-Zoeller. Up to duality, we have 
only 10 possible types. 
Next we prove that it is possible to classify all Hopf algebras 
of 8 of the 10 possible types.
 
\begin{theorem}\label{tabla}
(a) There is no Hopf algebra $H$ of dimension $p^{3}$ 
such that $H$ or $H^{*}$ is of one of the following
types:
$$(1; 1),\ (p, p; 1),\  
(p, p; p),\ (p, p; p^{2}),\ (p^{2}; 1).$$
(b) Let $H$ be a non-semisimple Hopf algebra of dimension $p^3$. 
\begin{enumerate}
\item If $H$ is of type 
$(p,p;p,p)$, then $H \simeq T(q) \otimes k[\ZZ/(p)]$.
\item If $H$ is of type $(p^{2}; p)$, then $H \simeq {\bf r}(q)$. 
\item If $H$ is of type $(p^{2}; p^{2})$, then either 
$H \simeq \widehat{T(q)}$ or $H \simeq \widetilde{T(q)}$. 
\end{enumerate}
\end{theorem}

\begin{proof}
It follows from Masuoka's classification that no semisimple Hopf 
algebra satisfies 
any of the conditions in $(a)$. Then, 
if there exists a Hopf algebra $H$ which satisfies one of the conditions
above, $H$ should be non-semisimple.
Type $(1;1)$ is not possible by Remark \ref{oddim}. 
For the other types, we have that  
$|G(H)| = p^2$
and the order of the antipode divides $4p^{2}$, by Radford's formula.
Hence, by Proposition \ref{schrad}, $H$ 
should be pointed. By inspection, the cases in $(a)$ are  
impossible and in case (1), $H$ should be isomorphic to  
$T(q)\otimes k[\ZZ/(p)]$, in case (2), $H$ should be isomorphic to ${\bf r}(q)$
and in case (3), $H$
should be isomorphic either to $\widehat{T(q)}$ or $\widetilde{T(q)}$.
\end{proof}

\begin{remark}
It follows from Theorem \ref{tabla} and Corollary \ref{simplicidad}  
that a Hopf algebra $H$ of dimension $p^{3}$ is isomorphic to a 
Hopf algebra of the list (a), \ldots, (k) or $H$ satisfies condition (I) or (II), 
{\it i.e.} $H$ is strange
and of type $(p;1)$ or $H$ is strange and of type $(p;p)$.
\end{remark}

\bigbreak
From Masuoka's classification it follows
that there is no semisimple Hopf algebra of dimension $p^{3}$  
of type $(p;1)$ or $(p;p)$. The Frobenuis-Lusztig kernels 
${\bf u}_{q}({\mathfrak {sl}}_{2})$, 
$q \in \GG_{p} \smallsetminus \{1\}$ are of type $(p; 1)$ and
the book algebras ${\bf h}(q,m)$, $q \in \GG_{p}\smallsetminus \{1\}$,
$m\in \ZZ /(p)\smallsetminus \{0\}$ are of type $(p;p)$. These are
the only non-semisimple pointed Hopf algebras of type
$(p;1)$ or $(p;p)$ (see \cite[Section 6]{AS1}).

\subsection{Hopf algebras of type ${\bf (p;p)}$}\label{pp}
In order to complete the classification of Hopf algebras of dimension $p^{3}$, we give
in this section some results on Hopf algebras of type $(p;p)$.

\smallbreak
The following theorem will be also useful in the quasitriangular case.
\begin{theorem}\label{evalua0}
Let $H$ be a finite-dimensional non-semisimple Hopf algebra such that
$G(H)$ is abelian, ${\mathcal S}^{4p} = \id$ and $<\alpha, x> =1$ for all
$x \in G(H)$, where $\alpha$ is the modular element of $H^{*}$. 
Assume further that there exists a surjective Hopf algebra map
$\pi: H \to L$, such that $\pi(x) =1$ for all $x \in G(H)$, where $L$ is
a semisimple Hopf algebra such that $\dim L = \dim H/p|G(H)|$. Then 
$\ord {\mathcal S} = 4p$.
\end{theorem}

\begin{proof}
Since $H$ is non-semisimple, the order of the antipode is bigger than 2 and divides
$4p$, by Radford's formula. Clearly it cannot be $p$ since the order is even.
Assume that ${\mathcal S}^{2p} = \id$.

\smallbreak
Let $e_{0},\ \ldots,\ e_{s}$ be the primitive central idempotents of $k[G(H)]$.
As in the proof of Proposition \ref{schrad}, 
we can write $1 = e_{0} + \cdots + e_{s}$ and hence 
$H = \bigoplus_{j=0}^{s} I_{j}$, where $I_j = He_{j}$. 
Since $G(H)$ is abelian, it follows that $\dim k[G(H)]e_{j} = 1$ for all 
$0\leq j\leq s$, 
and this implies that 
$\dim I_{j} = \dim H / |G(H)|$ for all $0\leq j\leq s$, 
since 
$$\dim I_{j} =\dim H\otimes_{k[G(H)]}k[G(H)]e_{j} = 
\dim_{k[G(H)]} H \dim k[G(H)]e_{j}.$$

\noindent It was also shown in the proof of Proposition \ref{schrad} 
that these spaces are invariant under the action of ${\mathcal S}^{2}$ and  
$\Tr (S^{2}|_{I_{j}}) = 0$. Let $q\in \GG_{p} \smallsetminus \{1\}$.
Since ${\mathcal S}^{2p} = \id$, 
by Lemma \ref{gral} $(a)$ follows  that
$I_{j} = \bigoplus_{m=0}^{p-1} I_{j,m}$ for all $0\leq j\leq s$, where
$I_{j,m} = \{ h \in I_{j} : {\mathcal S}^{2}(h) = q^{m}h \}$, and 
$p\dim I_{j,m} = \dim I_{j}$ 
for all $0\leq m \leq p-1$.
In particular, 
$$\dim I_{0,0} = \dim I_{0}/p =\dim H/p |G(H)| = \dim L,$$ 
\noindent and we can decompose $H$ as
$H = \bigoplus_{j,m = 0}^{p-1} I_{j,m}$.

\smallbreak
Since $\pi(x) = 1$ for all $x \in G(H)$, it follows that $\pi(e_{j}) = 0$ for all
$j \neq 0$, that is $I_{j} \subseteq \Ker \pi$ for all $j\neq 0$.
Moreover, $I_{0,m} \subseteq \Ker \pi$ for all $m\neq 0$, since for all 
$h \in I_{0,m}$, with $m\neq 0$, we have that 
$q^{m}\pi(h) = \pi({\mathcal S}^{2}(h))$ 
$= {\mathcal S}^{2}(\pi(h)) = \pi(h)$, because $L$ is semisimple. Therefore 
$\bigoplus_{(j,m)\neq (0,0)} I_{j,m} \subseteq \Ker \pi$, which implies that
\begin{equation}\label{ker}
\Ker \pi = \bigoplus_{(j,m)\neq (0,0)} I_{j,m}, 
\end{equation}

\noindent since both have the same dimension.
\smallbreak
Let $\Lambda \in H$ be a non-zero left integral, then $\Lambda \in I_{0,0}$.
Indeed, since $H = \bigoplus_{j=0}^{s} I_{j}$, there exist $h_{0},\ \ldots,\ h_{s}$ in $H$
such that $\Lambda = h_{0}e_{0} + \cdots + h_{s}e_{s}$. Hence 
$\Lambda = $ $<\alpha, e_{0}> \Lambda = \Lambda e_{0}  = h_{0}e_{0}$, since
$<\alpha, x> =1$ for all $x \in G(H)$ and $e_{0} = \frac{1}{|G(H)|} \sum_{x\in G(H)} x$.
Moreover, by \cite[Prop. 3, d)]{R3}, we know that 
${\mathcal S}^{2}(\Lambda) = <\alpha, g^{-1}>\Lambda$, 
where $g \in G(H)$ is the modular element of $H$. This implies that
$\Lambda \in I_{0,0}$,
since by assumption $<\alpha, g^{-1}> = 1$ and $\Lambda= h_{0}e_{0} \in I_{0}$.

\smallbreak
On the other hand, since $H$ is non-semisimple, we have that 
$<\eps, \Lambda>$ 
$= <\eps, \pi(\Lambda)>$ 
$= 0$ and this implies that $\pi(\Lambda)=0$, since
$\pi(\Lambda)$ is a left integral in $L$ and $L$ is semisimple.
Therefore, $\Lambda \in \Ker \pi \cap I_{0,0}$, implying by (\ref{ker})
that $\Lambda = 0$, which is a contradiction to our choice of $\Lambda$.
\end{proof}

\begin{remark}
Let $H,\ I_{0},\ \Lambda$ and $\alpha \in G(H^{*})$ be as in 
Theorem \ref{evalua0} and let 
$e_{0},\ \ldots,\ e_{s}$, be the primitive idempotents of $k[G(H)]$. Then
$\Lambda \in I_{0}$ if and only if $<\alpha, x> = 1$, for all $x \in G(H)$.
\begin{proof}

Suppose that $\Lambda \in I_{0}$. Then there exists $h \in H$ such that 
$\Lambda = h e_{0}$. In particular, for all $x \in G(H)$ we have that
$\Lambda x = <\alpha, x>\Lambda = h e_{0}x = h e_{0}= \Lambda$, and this implies
that $<\alpha, x> = 1$, for all $x \in G(H)$.

\smallbreak 
Conversely, suppose that $<\alpha, x> = 1$, for all $x \in G(H)$.
Since $H = He_{0}+ \cdots + He_{s}$, there exist $h_{0},\ \ldots,\ h_{s} \in H$ such
that $\Lambda = h_{0}e_{0}+ \cdots + h_{s}e_{s}$. 
Since $<\alpha, e_{0}>= 1$, it follows that 
$\Lambda = <\alpha, e_{0}>\Lambda $ $ = \Lambda e_{0}$ $= h_{0}e_{0}$, which implies that
$\Lambda \in I_{0}$. 
\end{proof}
\end{remark}

\begin{remark} The proof of the theorem above was 
inspired in some results of Ng;  
the spaces $I_{j,m},\ 0\leq j,m\leq p-1$ are the spaces 
$H^{w}_{0,m,j},\ w = q \in \GG_{p} \smallsetminus \{1\}$, 
defined in \cite[Section 3]{Ng} in the special case where
${\mathcal S}^{2p} = \id$.
\end{remark}

\begin{corollary}\label{evalua}
Let $H$ be a non-semisimple Hopf algebra of dimension $p^{3}$ 
and type $(p;p)$. 
\begin{enumerate}
\item Then the order of the antipode is $2p$ or $4p$.
\item If $<\beta, x> =1$, for all $\beta \in G(H^{*}),\ x\in G(H)$, 
then the order of the antipode is $4p$.
\end{enumerate}
\end{corollary}

\begin{proof}
(1) Since $H$ is non-semisimple, the order of the antipode
is bigger than 2 and divides $4p$, by Radford's formula. 
Since $\ord {\mathcal S}$ is even and $p$ is odd, 
it is necessarily $2p$ or $4p$. 
\smallbreak
(2) Consider the surjective Hopf algebra
map $H \xrightarrow{\pi} k^{G(H^{*})}$, defined by 
$<\pi(h), \beta> = <\beta, h>$, for all $h \in H$, $\beta \in G(H^{*})$.
Then $\pi(x) = 1$ for all $x\in G(H)$, since by assumption $<\beta, x> = 1$, 
for all $\beta \in G(H^{*}),\ x\in G(H)$. The claim follows directly from 
Theorem \ref{evalua0}, since
$k^{G(H^{*})}$ is semisimple, $\dim k^{G(H^{*})} = p = 
\dim H / p |G(H)|$ and ${\mathcal S}^{4p} = \id$.
\end{proof}

\bigbreak
Let $H$ be a Hopf algebra provided with a projection 
$H \xrightarrow{\pi} B$,
which admits a section of Hopf algebras 
$B \xrightarrow{\gamma} H$. Then 
$A = H^{co\ \pi}$ is a Hopf algebra in the category 
of Yetter-Drinfel'd modules 
over $B$ and $H$ is isomorphic to the smash product $A\# B$. 
In this case,
following the terminology of Majid, we say that 
$H$ is a {\it bosonization} of $B$. 
For references on the correspondence between Hopf 
algebras with a projection and 
Hopf algebras in the category of Yetter-Drinfel'd 
modules see \cite{Mj}, \cite{R5}.

\begin{proposition}\label{boson}
Let $H$ be a finite-dimensional non-semisimple Hopf algebra.
Assume that $G(H)$ is non-trivial and abelian and  
$G(H) \simeq G(H^{*})$. 
\begin{enumerate}
\item If $|G(H)|= p$, then $H$ is a bosonization of $k[G(H)]$ if and 
only if there exist $\beta \in G(H^{*})$ and $x\in G(H)$ such that 
$<\beta, x> \neq 1$.
\item If $\dim H = p|G(H)|$, then $H$ is a bosonization of $k[G(H)]$. 
\end{enumerate}
\end{proposition}

\begin{proof}
$(1)$ Suppose that $H$ is a bosonization of $k[G(H)]$. 
Then there exists a Hopf 
algebra projection $H \xrightarrow{\pi} k[G(H)]$
which admits a Hopf algebra section 
$k[G(H)] \xrightarrow{\gamma} H$, that
is $\pi \circ \gamma = \id$. 
Let $G = G(H)$ and denote by 
$\hat{G}$ the group of characters of
$G$.
We identify $\hat{G}$ with $\Alg (k[G],k) = G(k[G]^{*})$.  
Since $G \neq 1$, there exist $x \in G,\ \beta \in \hat{G}$ 
such that
$<\beta, x> \neq 1$. Consider the group homomorphism 
$\hat{G} \xrightarrow{\varphi} G(H^{*})$
given by
$\varphi(\chi)= \chi \circ \pi$, for all $\chi \in \hat{G}$.
Then $\beta \circ \pi \in G(H^{*})$, $\gamma(x) \in G$ and 
$<\beta \circ \pi, \gamma(x)>$ $= <\beta , \pi \circ \gamma(x)> =$ 
$ <\beta , x> \neq 1.$

\smallbreak 
Suppose now that $H$ is not a bosonization of $k[G]$. We show that 
for all $\beta \in G(H^{*}),\ x\in G$,
$<\beta, x> = 1$.
Denote by $k[G] \xrightarrow{\iota} H$ the inclusion of the group algebra in $H$.
Since $G(H^{*})$ is also non-trivial, we have a surjection of Hopf algebras
$H \xrightarrow{\pi} k^{G(H^{*})}$, given by $<\pi(h), \beta> = <\beta, h>$, 
for all $h \in H$, $\beta \in G(H^{*})$.

\smallbreak
\noindent{\bf{Claim 1.}} $\pi (x) = 1_{k^{G(H^{*})}}$, for all $x \in G$.
Since $G(H^{*})$ is abelian, we have that 
$k^{G(H^{*})} \simeq k[\widehat{G(H^{*})}] \simeq k[G]$ as Hopf algebras. Moreover, 
the composition of these isomorphisms with $\pi$ induce a 
Hopf algebra surjection $H \xrightarrow{\tau ' } k[G]$.
If there exists $x \in G$ such that 
$\pi(x) \neq 1_{k^{G(H^{*})}}$, then 
$\tau '(x) \neq 1$ and the restriction 
$\tau '\circ \iota$ of $\tau '$ 
to $k[G]$ defines an automorphism of $k[G]$, 
since $\vert G\vert = p$,
which is a prime number. Define now $\tau : H \to k[G]$ via
$\tau = (\tau '\circ \iota)^{-1} \circ \tau '$. 
Then $\tau \circ \iota = \id_{k[G]}$ and $H$
is a bosonization of $k[G]$, which is a contradiction to our assumption.

\smallbreak 
Therefore for all $\beta \in G(H^{*}),\ x \in G$ we have
$<\beta, x> = <\pi(x), \beta> = $
$< 1_{k^{G(H^{*})}}, \beta> =1$, 
which proves $(1)$.

\smallbreak
$(2)$ Consider the Hopf algebra surjection $H \xrightarrow{\tau '} k[G]$ defined in the
proof of Claim 1. 
\smallbreak
\noindent{\bf{Claim 2.}} $\tau'|_{G}$ defines a group automorphism of $G$.
Suppose on the contrary, that $\tau '|_{G}$ does not define an
automorphism of $G$. Then there exists $h \in G$ such that 
$h \neq 1$ and $\tau '(h) = 1$; in particular, $h \in H^{co\ \tau'}$.

\smallbreak
On the other hand, 
$\dim H^{co\ \tau'}
= p$, since 
$\dim k^{G(H^{*})} = |\widehat{G(H^{*})}| = |G|$, and $\dim H = p|G| = $
$\dim H^{co\ \tau'}\dim k^{G(H^{*})}$, by Nichols-Zoeller.
Since $p$ is a prime number, it follows that $\ord h = p$ and
$k<h> = H^{co\ \tau '}$. In particular, one has the exact sequence
of finite-dimensional Hopf algebras
$ 1\to k<h> \xrightarrow{\iota} H \xrightarrow{\tau '} 
k[G] \to 1, $ 
implying as it is well-known that $H$ is semisimple, 
which is a contradiction to our assumption.

\smallbreak 
Therefore the automorphism $\tau '|_{G}$ defines a Hopf algebra automorphism
in $k[G]$, and $H$ is a bosonization of $k[G]$.   
\end{proof}

\begin{remark}\label{bosones}
If we examine the Hopf algebras in the list (a), \ldots, (k) in the
introduction, we see
that the cases (d), (e), (f) and (i) are also bosonizations.
In the case (d) it is clear that the product Hopf algebra is the 
bosonization of $k[\ZZ /(p)]$ and by the proposition above it is
also a bosonization of $k[\ZZ /(p) \times \ZZ /(p)]$. The cases (e)
and (f) are bosonizations of $k[\ZZ /(p^{2})]$ and the book algebras
${\bf h} (q,m)$ in case (i) are bosonizations of $k[\ZZ /(p)]$.
In the case of the book algebras, Andruskiewitsch and Schneider
proved in \cite{AS1} that these algebras are also 
bosonizations of a Taft algebra, {\it i.e.} they are isomorphic
to a smash product $R\# T(q)$. Moreover, they also proved  
that the algebras $R$ exhaust the list of non-semisimple
Hopf algebras of order $p$ in the Yetter-Drinfel'd category over $T(q)$.      
\end{remark}

\bigbreak
We now prove that a Hopf algebra $H$ of dimension $p^{3}$ is a bosonization
of $k[\ZZ /(p)]$ in the following two cases. Although some properties of
$H$ are known, we cannot determine its structure, since 
the classification 
of Hopf algebras of dimension $p^{3}$ which are bosonizations of $k[\ZZ /(p)]$, 
that is, of braided Hopf algebras of dimension $p^{2}$ in the category of
Yetter-Drinfel'd modules over $k[\ZZ /(p)]$, is not known. 
 
\begin{corollary}\label{boson1}
Let $H$ be a non-semisimple Hopf algebra of dimension $p^{3}$ and type
$(p;p)$. If 
${\mathcal S}^{2p} = \id$ then $H$ is a bosonization of $k[\ZZ /(p)]$.
\end{corollary}

\begin{proof}
Follows from Corollary \ref{evalua} (2) and Proposition \ref{boson} (1).
\end{proof}

\begin{corollary}\label{usable}
Let $H$ be a Hopf algebra of dimension $p^{3}$ and type
$(p;p)$ such that $H$ contains a non-trivial skew primitive
element. Then
$H$ is a bosonization of $k[\ZZ /(p)]$.
\end{corollary}

\begin{proof}
Suppose on the contrary that $H$ is not 
a bosonization of $k[\ZZ /(p)]$.
By \cite[Prop. 1.8]{AN}, $H$ contains a
Hopf subalgebra $B$ which is isomorphic to a 
Taft algebra of dimension $p^{2}$.
Let  $x,\ h \in H$, be the generators of $B$,
such that $1 \neq h\in G(H)$ and $x \in P_{1,h}$. 

\smallbreak
Consider now the Hopf algebra surjection
$\pi: H \xrightarrow{} k^{G(H^{*})}$ given by\\ $<\pi(t),\beta>$
$=<\beta,t>,$
for all $t \in H,\ \beta \in G(H^{*})$. 
If $\pi$ is not trivial on $G(H)$, then $H$ is a bosonization of  
$k[\ZZ /(p)]$, by Proposition \ref{boson} (1).
If $\pi$ is trivial on $G(H)$, then 
$t \in H^{co\ \pi}$, for all $t\in G(H)$.
Since $\pi$ is a Hopf algebra map and 
$\Delta (x) = x\otimes 1 + h\otimes x$, 
we see that $\pi(x)$ is a primitive element in $k^{G(H^{*})}$. 
Since $G(H^{*})$ has finite order,  
$\pi(x) = 0$ and therefore $(\id\otimes \pi)\Delta(x) = x\otimes 1$, that is
$x \in H^{co\ \pi}$.
Since $B$ is generated by $x$ and
$h \in G(H)$, we have that $B \subseteq H^{co\ \pi}$ and hence
$B= H^{co\ \pi}$, because both have dimension $p^{2}$.
Therefore $H$ fits into the extension of Hopf algebras 
$$1\xrightarrow{} B \xrightarrow{\iota} H \xrightarrow{\pi}k^{G(H^{*})} \to 1.$$

\noindent As $G(H^{*}) \simeq \ZZ /(p)$, $H^{*}$ is pointed by Theorem \ref{dim1}.
Since the only non-semisimple pointed Hopf algebras of dimension
$p^{3}$ of type $(p;p)$ are the book algebras, it follows that 
$H^{*} \simeq {\bf h}(q,-m)$, where
$q \in \GG_{p} \smallsetminus \{1\},$ and $m \in \ZZ /(p)\smallsetminus \{0\}$. 
Hence  
$H \simeq {\bf h}(q,m)$ and therefore pointed, since
for all $q \in \GG_{p} \smallsetminus \{1\},\ m \in \ZZ /(p)\smallsetminus \{0\}$ 
we have that ${\bf h}(q,-m)^{*} \simeq {\bf h}(q,m)$
(see \cite[Section 6]{AS1}). Hence $H$ is a bosonization of $k[\ZZ /(p)]$, by
Remark \ref{bosones}. 
\end{proof}

Finally we note a very special result on Hopf algebras of dimension $p^{3}$ following
\cite[Prop. 1.3]{Na2} and Theorem \ref{dim1}.

\begin{proposition}
Let $H$ be a non-semisimple Hopf algebra of dimension $p^{3}$ and 
assume that $H$ contains a simple subcoalgebra $C$
of dimension 4 such that ${\mathcal S}(C) = C$. Then $H^{*}$ is
pointed. In particular, $H$ cannot be of type $(p;p)$.
\end{proposition}

\begin{proof}
Let $B$ be the algebra generated by $C$; clearly $B$ is a Hopf subalgebra 
of $H$ and it follows that $\dim B \vert p^{3}$, by Nichols-Zoeller. 
Since $C \subseteq B \subseteq H$, $B$ is non-semisimple and non-pointed. 
Then necessarily $B = H$, 
since by \cite[Thm. 2]{Z} and \cite[Thm. 5.5]{Ng}, the only non-semisimple 
Hopf algebras whose dimension is a power of $p$ with exponent 
less than 3 are the Taft algebras, which are pointed.
Hence $H$ is generated
as an algebra by a simple coalgebra of dimension 4 which is stable by the antipode.
By \cite[Prop. 1.3]{Na2}, $H$ fits into an extension 
$1\to k^{G} \to H \to A \to 1$, where $G$ is a finite group and $A^{*}$ is a 
pointed non-semisimple Hopf algebra. 
\smallbreak 
Since $H$ is not semisimple, it follows that 
$|G(H)| = 1$ and $H^{*}$ is pointed, or $|G(H)|= p$ and $H$ is pointed 
by Theorem \ref{dim1}, which is impossible by assumption. Moreover, 
if $H^{*}$ is pointed and of type $(p;p)$, then $H^{*}$ is isomorphic to book algebra 
${\bf h}(q,m)$, 
for some $q \in \GG_{p}\smallsetminus \{1\},\ m \in \ZZ /(p)\smallsetminus \{0\}$. 
Hence $H$ is also pointed
and it cannot contain a simple subcoalgebra of dimension 4. 
\end{proof}

\section{Quasitriangular Hopf algebras of dimension $p^3$}\label{qt}
Let $H$ be a finite-dimensional Hopf algebra and let 
$R \in H\otimes H$. As usual, we use for $R$ the symbolic notation 
$R = R^{(1)} \otimes R^{(2)} $. Define a linear map 
$f_{R}: H^{*} \to H$ by
$f_{R}(\beta) =  <\beta,R^{(1)}>R^{(2)}$, for $\beta \in H^{*}$. The pair $(H,\ R)$ 
is said to be a quasitriangular Hopf algebra \cite{Dr}
if the following axioms hold:

\smallbreak
(QT.1) $\Delta^{\cop}(h)R = R \Delta(h)$, $\forall\ h \in H$,
\smallbreak
(QT.2) $(\Delta \otimes \id)(R) = R_{13}R_{23}$,
\smallbreak
(QT.3) $(\eps \otimes \id)(R) = 1$,
\smallbreak
(QT.4) $(\id \otimes \Delta)(R) = R_{13}R_{12}$,
\smallbreak
(QT.5) $(\id \otimes \eps)(R) = 1$;

\smallbreak
\noindent or equivalently if $f_{R}: H^{*cop} \to H$ is a bialgebra map and 
(QT.1) is satisfied.

\smallbreak
We have used the notation $R_{12}$ to indicate the element 
$R \otimes 1 \in H^{\otimes 3}$, similarly for $R_{13}$ and $R_{23}$. 
Note that $(H^{\cop}, R_{21})$ and $(H^{\op},
R_{21})$ are also quasitriangular, where 
$R_{21} : = R^{(2)}\otimes R^{(1)}$.

\smallbreak
We refer to a pair $(H,\ R)$ which satisfies the five axioms 
above as a quasitriangular Hopf algebra or simply saying that 
$H$ admits a quasitriangular structure.
See for example \cite{R2} for a fuller discussion and references.

\smallbreak
We define a morphism $f: (H,\ R) \to (H',\ R')$ 
of quasitriangular Hopf algebras over $k$
to be a Hopf algebra map $f: H \to H'$ such that $R' = (f \otimes f)(R).$

\smallbreak
Let $\widetilde{R} = R_{21}$. There is another Hopf algebra map 
$f_{\widetilde{R}}: H^{*} \to H^{op}$,
given by $f_{\widetilde{R}}(\beta) = <\beta,R^{(2)}> R^{(1)}$, for all $\beta \in H^{*}$. 
With the usual identification
of vector spaces of $H$ and $H^{**}$, the maps 
$f_{\widetilde{R}}$ and $f_{R}$ are related by
the equation $f_{\widetilde{R}} = f_{R}^{*}$.

\smallbreak
\begin{remark}\label{KL}
Let $L$ and $K$ denote, respectively, the images of $f_R$ and 
$f_{\widetilde{R}}$.
Then $L$ and $K$ are Hopf subalgebras of $H$ of dimension $n > 1$, 
unless $H$ is
cocommutative and $R = 1 \otimes 1$; this dimension is called the 
rank of $R$.
By \cite[Prop. 2]{R2}, we have  $L \simeq {K^*}^{\cop}$.

\smallbreak
Let $H_{R}$ be the Hopf subalgebra of $H$ generated by $L$ and $K$. 
If $B$ is a Hopf
subalgebra of $H$ such that $R \in B\otimes B$, then $H_{R} \subseteq B$. Hence,
we say that $(H,\ R)$ is a {\it{minimal}} quasitriangular Hopf algebra if 
$H = H_{R}$.
It is shown in \cite[Thm. 1]{R2} that $H_R = LK = KL$.
If $L$ is semisimple, then $K$ is semisimple
and therefore $H_R$ is semisimple.
Minimal quasitriangular Hopf algebras were first 
introduced and studied in \cite{R2}.
\end{remark}

\bigbreak

We recall some fundamental properties of finite-dimensional 
quasitriangular Hopf algebras (see for example 
\cite{Dr}, \cite{Mo}). 
First $R$ is invertible with inverse
$R^{-1} = ({\mathcal  S}\otimes I)(R) = (I\otimes {\mathcal  S}^{-1})(R)$, 
and $R = ({\mathcal  S}\otimes {\mathcal  S})(R)$.
Set $u = {\mathcal  S}(R^{(2)}) R^{(1)}$, then $u$ 
is also invertible where
\begin{align*}
u^{-1} = R^{(2)}{\mathcal  S}^{2}(R^{(1)}),
\ \Delta (u) & = (u\otimes u)(\widetilde{R}R)^{-1}= 
(\widetilde{R}R)^{-1}(u\otimes u),\\
\eps (u)  =  1
\mbox{ and }
{\mathcal  S}^{2}(h)  = uhu^{-1} & =  
({\mathcal  S}(u))^{-1}h{\mathcal  S}(u), \mbox{ for all }h \in H. 
\end{align*}

Consequently, $u{\mathcal  S}(u)$ 
is a central element of $H$.
Since ${\mathcal  S}^{2}(h) = h$, for all $h \in G(H)$, it follows that 
$u$ commutes with the group-like elements of $H$. 
The element $u \in H$ is called the {\it{Drinfel'd
element}} of $H$.

\smallbreak
We say that $v \in H$ is a {\it ribbon element} of $(H,\ R)$ if the 
following conditions are satisfied:

\smallbreak
(R.1) $v^{2} = u{\mathcal S}(u),$
\smallbreak
(R.2) ${\mathcal S}(v) =v,$
\smallbreak
(R.3) $\eps (v) = 1,$
\smallbreak
(R.4) $\Delta (v) = (\widetilde{R}R)^{-1}(v\otimes v)$ and
\smallbreak
(R.5) $vh = hv,$ for all $h \in H$.

\smallbreak
If $H$ contains a ribbon element, then the triple $(H,\ R,\ v)$ or simply $H$ 
is called a {\it ribbon Hopf algebra} 
(see \cite[XIV. 6]{K}, \cite{KR}, \cite[Section 2.2]{R4}). 

\smallbreak
The following theorem is due to Natale and it will be crucial to prove
some results in the case where the dimension of the Hopf algebra is $p^{3}$. 

\begin{theorem}\cite[Thm. 2.3]{Na1}\label{Natale}
Let $H$ be a Hopf algebra of dimension $pq$ over $k$, where $p$ and $q$ 
are odd primes which are not necessarily distinct. Assume that $H$ admits a 
quasitriangular structure. Then $H$
is semisimple and isomorphic to a group algebra $k[F]$, 
where $F$ is a group of order $pq$.
\end{theorem}

\begin{remark}\label{taftnqt}
The theorem above implies the known fact that the
Taft algebra $T(q)$ does not admit any 
quasitriangular structure if $\dim T(q) = p^{2}$, with $p$ odd prime.
\end{remark}

\bigbreak
The following result is due to Gelaki.

\begin{theorem}\label{Gel} 
Let $(H,\ R)$ be a finite-dimensional 
quasitriangular Hopf algebra with antipode
${\mathcal  S}$ over a field $k$ of characteristic 0.

\smallbreak 
\begin{enumerate}
\item[(a)] \cite[Thm. 3.3]{Ge1} If the Drinfel'd
element $u$ of $H$ acts as a scalar in any irreducible 
representation of $H$ (e.g. when $H^{*}$ is pointed), 
then $u = {\mathcal  S}(u)$ and in particular
${\mathcal  S}^{4}=\id$.
\item[(b)] \cite[Thm. 1.3.5]{Ge2}
If $H_{R}$ is semisimple, then $u = {\mathcal S}(u)$ and 
${\mathcal S}^4 = \id$.
\item[(c)] If $H^{*}$ is pointed then either $H$ 
is semisimple or $\dim H$ is even.
\end{enumerate}
\end{theorem}

\begin{proof}
(c) follows from (a) and Remark \ref{oddim}.
\end{proof}

\bigbreak
We can now prove our first assertion.

\begin{proposition}\label{uq}
Among the Hopf algebras in the list (a), \ldots , (k) in the
introduction, only
the group algebras and the 
Frobenius-Lusztig kernels ${\bf u}_{q}({\mathfrak {sl}}_{2})$, where 
$q \in \GG_{p} \smallsetminus \{1\}$,  admit a quasitriangular structure.
\end{proposition}

\begin{proof}
It is well-known that the group algebras and the
Frobenius-Lusztig kernels are quasitriangular (see \cite[IX.7]{K}). 
We show next
that the other Hopf algebras in the list 
cannot admit a quasitriangular structure. 

\smallbreak
Let $q\in \GG_{p}\smallsetminus \{1\}$.
The Hopf algebras in cases (d), (f) and (g)
admit no quasitriangular structure, since they have a
Hopf algebra surjection to the Taft algebra $T(q)$ and by Remark \ref{taftnqt},
$T(q)$ is not quasitriangular.
(see \cite[Section 1]{AS2}). 

\smallbreak
Let $H$ be one of the Hopf algebras $\widetilde{T(q)}$,
${\bf h} (q, m),\ m\in \ZZ/(p)\smallsetminus \{1\}$, 
${\bf u}_{q}({\mathfrak {sl}}_{2})^{*}$, ${\bf r}(q)^{*}$, {\it i.e.} 
$H$ is one of the cases (e), (i), (j), (k) of the list. 
Then $H^{*}$ is pointed and $H$ is not semisimple of odd dimension. 
Hence by Theorem \ref{Gel} (a) and Remark \ref{oddim},
$H$ cannot admit a quasitriangular structure. 

\smallbreak
Let $G$ be a finite group. If $H= k^{G}$ admits a quasitriangular structure, 
then $G$ must be abelian by (QT1), and 
$H$ is isomorphic to a group algebra.

\smallbreak
Finally, the semisimple Hopf algebras of dimension $p^3$ in (c) 
are not quasitriangular by \cite[Thm. 1]{Mk4}.
\end{proof}

\bigbreak
\begin{remark}
The case of $\widetilde{T(q)}$ of Proposition \ref{uq} also follows from
\cite[Section 5]{R4}. There, Radford defines Hopf algebras which depends 
on certain parameters. He proved
that these algebras admit a quasitriangular structure if and only 
if these parameters satisfy specific
relations. One can see that the algebras  are of 
this type and
the relations needed to have a quasitriangular structure do not hold.
\end{remark}

\bigbreak
In the following, we give a partial description of the 
quasitriangular Hopf algebras of dimension $p^{3}$.

\smallbreak
Let $D(H)$ be the Drinfel'd double of $H$ (see \cite{Mo} or \cite{R2} 
for its definition and properties). We identify 
$D(H) = H^{* cop} \otimes H$ as vector spaces and for 
$\beta \in H^{*}$ and $h \in H$, the element 
$\beta \otimes h \in D(H)$ is denoted by $\beta \# h$.
We may also identify 
$D(H)^{*} = H^{op} \otimes H^{*}$, and for 
$h \in H^{op}$ and $\beta \in H^{*}$, the element 
$h \otimes \beta  \in D(H)^{*}$ is denoted by $h\# \beta$ if
no confusion arrives.

\smallbreak
For every finite-dimensional quasitriangular Hopf algebra $(H,\ R)$, 
there is a Hopf algebra surjection 
$$D(H) \xrightarrow{F} H,\quad F(\beta \# h) =  
<\beta, R^{(1)}> R^{(2)}h,$$ 

\noindent which induces by duality 
an inclusion of Hopf algebras $H^{*} \xrightarrow{F^{*}} D(H)^{*}$. 
Moreover, by \cite[Prop. 10]{R2} all the
group-like elements of $D(H)^{*}$ have the form
$x\# \beta$, for some $x \in G(H),\ \beta \in G(H^{*})$,
and there is a central extension of Hopf algebras
\begin{equation}\label{natext}
1\to k[G(D(H)^{*})] \xrightarrow{\iota} D(H) \xrightarrow{\pi}
A \to 1,
\end{equation}

\noindent where $\iota (x\# \beta) = \beta\# x$, 
for all $x\# \beta \in G(D(H)^{*})$
and $A$ is the Hopf algebra given by the quotient 
$D(H) / D(H) k[G(D(H)^{*})]^{+}$.

\smallbreak
We need the following lemma due to Natale.

\begin{lemma}\label{sonia}\cite[Lemma 3.2]{Na1}
Let $1 \to A \xrightarrow{\iota} H \xrightarrow{\pi} B \to 1,$ be an extension
of finite-dimensional Hopf algebras. Let also $L \subseteq H$ be a Hopf
subalgebra. If $L$ is simple then either $L \subseteq \iota(A)$ or $L \cap \iota(A) = k1$.
In the last case, the restriction $\pi|_{L} : L \to B$ is injective.
\end{lemma}

\smallbreak
\begin{lemma}\label{uno}
Let $H$ be a quasitriangular Hopf algebra of dimension $p^{3}$ 
such that the Hopf algebra map $f_R$ defined above is an isomorphism. 
Assume further that $H$ is simple as a Hopf algebra.
Then $G(H)$ has order $p$ and $<\alpha, g> = 1$, where $\alpha$ and $g$ are 
the modular elements of $H^{*}$ and $H$, respectively. 
\end{lemma}

\begin{proof}
$H$ cannot be semisimple, since otherwise $H$ would admit a
non-trivial central group-like element by \cite[Thm. 1]{Mk3} 
and this would contradict
our assumption on the simplicity of $H$.
Moreover, $H$ cannot be unimodular, 
since otherwise $H^{*}$
would be also unimodular and by Radford's formula, 
${\mathcal S}^{4} = \id$,
which implies by Remark \ref{oddim} that $H$ is semisimple. 
Since $f_{R}$ is an isomorphism, 
$H^{*}$ is also non-semisimple, simple as a Hopf algebra and 
$G(H)\simeq G(H^{*})$. Hence, by Nichols-Zoeller
we have to deal only with the cases where 
$|G(H)| = |G(H^{*})| = p$ or $|G(H)| = |G(H^{*})| = p^{2}$.
But the order of $G(H)$ cannot be $p^{2}$, since otherwise $H$ 
would be pointed by Proposition \ref{schrad} and hence
isomorphic to a Frobenius-Lusztig kernel 
${\bf u}_{q}({\mathfrak {sl}}_{2})$, 
for some $q\in \GG_{p} \smallsetminus \{1\}$, by Proposition \ref{uq}, 
which is a contradiction, since 
$\vert G({\bf u}_{q}({\mathfrak {sl}}_{2}))\vert= p$. 
Therefore, the only possibility is $|G(H)| = |G(H^{*})| = p$.

\smallbreak
By \cite[Cor. 2.10, 1)]{Ge2}, we know that 
$f_{\widetilde{R}} (\alpha) = g^{-1}$ and by \cite[Prop.3]{R4}, 
$f_{\widetilde{R}R} = f_{\widetilde{R}} * f_{R}$ and
$f_{\widetilde{R}R} (\alpha) = 1$. This implies that necessarily 
$f_{R} (\alpha )= g$ and hence the order of $g$ and $\alpha$ 
must be equal.

\smallbreak
Consider now the Hopf algebra surjection
$F: D(H) \to H$ and the Hopf algebra inclusion $F^{*}: H^{*} \to D(H)^{*}$,
defined above. Since $G(H^{*}) \neq 1$, we have that 
$G = G(D(H)^{*}) \neq 1$, and dualizing extension (\ref{natext}), 
we get another extension of Hopf algebras given by
\begin{equation}
1\to A^{*} \xrightarrow{\pi^{*}} D(H)^{*} \xrightarrow{\iota^{*}}
k^{G} \to 1.
\end{equation}

Define $L = F^{*}(H^{*}) \subseteq D(H^{*})$. 
Since $H^{*}$ is simple as a Hopf algebra, 
by Lemma \ref{sonia}
we have that $L \subseteq \pi^{*}(A^{*})$ or $L\cap \pi^{*}(A^{*}) = k1$.
But the last case cannot occur, since otherwise the restriction
$\iota^{*}|_{L}: L \to k^{G}$ would be injective, implying
that $H^{*}$ is semisimple. Hence $L \subseteq \pi^{*}(A^{*})$. 

\smallbreak
Then there are $\beta \in G(H^{*}),\ \beta \neq \eps$ and 
$x \in G(H) \smallsetminus \{1\}$ such 
that $x \# \beta \in G(\pi^{*}(A^{*})) \subseteq G$. Moreover, since the image of
$G$ in $D(H)$ is central and $F$ is a Hopf algebra surjection, 
it follows that $F(\beta\# x) \in G(H) \cap Z(H)$ 
and by simplicity of $H$ we have that necessarily $F(\beta\# x) = 1$.

\smallbreak
Since both modular elements are not trivial, 
we have that $G(H^{*}) =$ $<\alpha>$ and $G(H) = <g>$ and hence 
$|G(D(H))| = p^{2}$.
Moreover, $|G(H^{*})|= |G(A^{*})|= p$, since otherwise 
$|G(A^{*})|= |G(D(H)^{*})|= |G(D(H))| = p^{2}$ and 
this would imply that $H$ has a central non-trivial group-like element, 
which contradicts our assumption on the simplicity of $H$.  
Hence the group-like element $\beta\#x$ generates $G(\pi^{*}(A^{*}))$, and
$\beta = \alpha^{j}$, $x=g^{i}$ for some
$1\leq i,j \leq p-1$. In particular, 
$$1 = F(\alpha^{j}\# g^{i}) = <\alpha^{j}, R^{(1)}> R^{(2)}g^{i} = 
f_{R} (\alpha^{j})g^{i} =  
f_{R} (\alpha)^{j}g^{i} = g^{j}g^{i}.$$

\noindent Then we have that $j \equiv -i \mod {p}$ and
$\pi^{*}(G(A^{*})) = <g \# \alpha^{-1}>$.   
Moreover, by \cite[Cor. 2.3.2]{Na2},
it follows that $<\alpha^{-1},g>^{2} = 1$, and this implies that 
$1 = <\alpha^{-1},g> =<\alpha,g>^{-1}$, since $|G(H)| = |G(H^{*})| = p$ and 
$p$ is odd.
\end{proof}
 
\bigbreak
We prove next one of our main results. 

\begin{theorem}\label{principal}
Let $H$ be a quasitriangular Hopf algebra of dimension $p^{3}$. Then
\begin{enumerate}
\item[$(i)$] $H$ is a group algebra, or
\item[$(ii)$] $H$ is a isomorphic to  
${\bf u}_{q}({\mathfrak {sl}}_{2})$, for some $q\in \GG_{p} \smallsetminus \{1\}$,
or
\item[$(iii)$] $H$ is a strange Hopf algebra of type $(p;p)$ 
and the map $f_{R}$ is an isomorphism. Moreover,
$H$ and $H^{*}$ are minimal quasitriangular, $1 = <\beta, x>$, for all
$\beta \in G(H^{*})$, $x\in G(H)$, and $\ord {\mathcal S} = 4p$.  
\end{enumerate}
\end{theorem}

\begin{proof}
If $H$ is semisimple, the claim follows by \cite[Thm. 1]{Mk4}, 
since $H$ must be isomorphic to a group algebra.

\smallbreak
Assume now that $H$ is non-semisimple and let $H_{R}$ be 
the minimal quasitriangular Hopf subalgebra
of $H$. 
Recall that $H_{R} = KL = LK,$ where 
$K$ = Im ${f_{R}}$ and $L$= Im $f_{\widetilde{R}}$. 
By Theorem \ref{Gel} (b) and Remark \ref{oddim}, 
$H_{R}$ is necessarily non-semisimple. 
Since the only non-semisimple Hopf algebra 
whose dimension is a power of $p$ with exponent 
less than 3 are the Taft algebras, 
by \cite[Thm. 2]{Z} and \cite[Thm. 5.5]{Ng}, and 
the Taft algebras are not quasitriangular by Remark \ref{taftnqt}, we
conclude that $\dim H_{R} = p^{3}$.

\smallbreak
Therefore the only possible case is when 
$H_{R} = H$ and $(H,\ R)$ is a minimal quasitriangular 
Hopf algebra. Then
by \cite[Cor. 3]{R2}, it follows that $\dim H \vert (\dim K)^{2}$;
hence the dimension of $K$ can only be $p^{2}$ or $p^{3}$.

\smallbreak
Suppose that the dimension of $K$ is $p^{2}$. 
Since $H$ is not semisimple, $K$ is not 
semisimple by Remark \ref{KL}. Moreover,  
by \cite[Thm. 5.5]{Ng} again, $K$ must be
isomorphic to a Taft algebra $T(q)$, where 
$q \in \GG_{p} \smallsetminus \{1\}$.
Since $L \simeq {K^*}^{\cop}$, $L$ must be also isomorphic to a 
Taft algebra and
by Remark \ref{taft}, $L \simeq T(q^{-1})$.

\smallbreak
It is clear that $G(K) \subseteq G(H)$ and $G(L) \subseteq G(H)$, 
where the order of $G(H)$
is $p$ or $p^{2}$. Since $H$ is a product of two Taft algebras, 
we have that ${\mathcal  S}^{4p} = \id$.
If the order of $G(H)$ is $p^{2}$, then by Proposition
\ref{schrad}, $H$ must be pointed. Then, by \cite[Thm. 0.1]{AS2} 
and Proposition \ref{uq},
$H$ is isomorphic to a Frobenius-Lusztig kernel 
${\bf u}_{q}({\mathfrak {sl}}_{2})$; but this cannot occur, since
$\vert G({\bf u}_{q}({\mathfrak {sl}}_{2}))\vert = p$.

\smallbreak
Therefore $\vert G(H)\vert = p$, and consequently $G(H) = G(K) = G(L)$.
Denote by $g \in G(H)$, and $x \in K$ the generators of $K$, and $g' \in G(H)$ 
and
$y \in L$ the generators of $L$; they are subject to similar
relations as in Remark \ref{taft}, since both are isomorphic to Taft algebras, 
but for different roots of unity. Moreover, 
$g'$ must be a power of $g$, since $\vert G(H)\vert = p$. 

\smallbreak
Hence $H$ is generated as an algebra by $g,\ x$ and $y$, which implies
that $H$ is pointed by \cite[Lemma 5.5.1]{Mo}.
Therefore, by \cite[Thm. 0.1]{AS2} and Proposition \ref{uq},
$H$ is isomorphic to ${\bf u}_{q}({\mathfrak {sl}}_{2})$,
for some $q\in \GG_{p} \smallsetminus \{1\}$.

\smallbreak
Assume now that $\dim K = p^{3}$, then the map 
$f_{R}: H^{*cop} \to H$ is an isomorphism of Hopf algebras. 
Moreover, $H$ must be non-pointed and hence simple as a 
Hopf algebra by Corollary \ref{simplicidad}, since otherwise
$H$ would be isomorphic to a Frobenius-Lusztig kernel 
${\bf u}_{q}({\mathfrak {sl}}_{2})$, 
which is a contradiction, since ${\bf u}_{q}({\mathfrak {sl}}_{2})^{*}$ is 
not pointed nor quasitriangular.
Then by Lemma \ref{uno}, 
we get that $|G(H)| = p$ and $<\alpha, g> =1$, where 
$\alpha $ and $g$ are the modular elements of $H^{*}$ and $H$ respectively.
It was shown in the proof of Lemma \ref{uno} that $f_{R}(\alpha) = g$; 
then it follows
by Remark \ref{oddim}, that $H$ and $H^{*}$ are both not unimodular.  
Hence, $<\beta, x> = 1$, for all $\beta \in G(H^{*})$, $x\in G(H)$, 
which implies by Corollary \ref{evalua} (2) that $\ord {\mathcal S} = 4p$.
Therefore, $H$ is a strange Hopf algebra of type $(p;p)$ which satisfies all 
conditions in $(iii)$.
\end{proof}

It is well-known that the group algebras and the 
Frobenius-Lusztig kernels are ribbon Hopf algebras (see \cite{K}). Next we prove
that there are no other ribbon Hopf algebras of dimension $p^{3}$. 

\begin{corollary}\label{ribbon}
Let $H$ be a ribbon Hopf algebra of dimension $p^{3}$. Then $H$ is a
group algebra or $H$ is isomorphic to 
${\bf u}_{q}({\mathfrak {sl}}_{2})$,
for some $q\in \GG_{p} \smallsetminus \{1\}$. 
\end{corollary}

\begin{proof}
Suppose on the contrary that $H$ is a ribbon Hopf algebra
of dimension $p^{3}$ which is not a group algebra and $H$ is not isomorphic to
a Frobenius-Lusztig kernel. Then by the preceding theorem, $H$ is of type
$(p;p)$ and $\ord {\mathcal S}= 4p$. But this cannot occur, 
since by \cite[Thm. 2]{KR}, the square of the antipode must have odd 
order. 
\end{proof}

As another application of Theorem \ref{principal} we classify quasitriangular
Hopf algebras of dimension 27 using some results from \cite{AN} and \cite{BD}.
As in \cite{BD}, we denote by $M^{c}(n,k)$ 
the simple matrix coalgebras contained in
the coradical.

\smallbreak
Let $H$ be a finite-dimensional Hopf algebra and
let $H_{0}$ be the coradical of $H$.
Then we have that $H_{0} \simeq \bigoplus_{\tau \in \hat{H}} H_{\tau}$, where
$H_{\tau}$ is a simple subcoalgebra of dimension $d_{\tau}^{2},\ d_{\tau} \in \ZZ$, and
$\hat{H}$ is the set of isomorphism types of simple left $H$-comodules. Define
$$H_{0,d} = \bigoplus_{\tau \in \hat{H}:\ d_{\tau} = d} H_{\tau}.$$
\noindent For instance $H_{0,1} = k[G(H)]$ and $H_{0,2}$ is the sum of all
4-dimensional simple subcoalgebras of $H$.
By \cite[Lemma 2.1 (i)]{AN}, the order of $G(H)$ divides the dimension of $H_{0,d}$
for all $d\geq 1$.  

\begin{theorem}
Let $H$ be a quasitriangular Hopf algebra of dimension 27. Then $H$ 
is a group algebra or $H$ is isomorphic to  
${\bf u}_{q}({\mathfrak {sl}}_{2})$,
for some $q\in \GG_{3} \smallsetminus \{1\}$. 
Precisely, $H$ is isomorphic to
one and only one Hopf algebra of the following list, where 
$q \in \GG_{3} \smallsetminus \{1\}$,
\begin{enumerate}
\item[(a)] $k[\ZZ / (27)]$
\item[(b)] $k[\ZZ /(9)\times \ZZ /(3)]$
\item[(c)] $k[\ZZ /(3) \times \ZZ /(3) \times \ZZ /(3)]$
\item[(d)] ${\bf u}_{q}({\mathfrak {sl}}_{2}) := k<g,\ x,\ y |\
gxg^{-1} = q^{2}x,\ gyg^{-1} = q^{-2}y,\ g^{3} = 1,\ x^{3} = 0,\
y^{3} = 0,\ xy - yx = g - g^{2}>.$ 
\end{enumerate}
\end{theorem}

\begin{proof}
By Theorem \ref{principal}, we have to show that there is no 
quasitriangular Hopf algebra  
$H$ of dimension 27 which satisfies $(iii)$.

\smallbreak
Suppose that such a Hopf algebra exists.
Since $H$ is not semisimple, $H$ is also not cosemisimple. 
Moreover, $H$ has no non-trivial skew primitive elements,
by Corollary \ref{usable} and Proposition \ref{boson} (1).

\smallbreak
Suppose now that 
$$H_{0} = k[G(H)] \oplus M^{c}(n_{1},k) \oplus \cdots \oplus M^{c}(n_{t},k),$$ 

\noindent where $2\leq n_{1} \leq \cdots \leq n_{t}\leq 3$, since 
$3 \vert \dim H_{0,d}$, for all $d\geq 1$. Moreover, since $H$ is not cosemisimple, 
we have the following possibilities
for $H_{0}$:

\begin{enumerate}
\item $H_{0} = k[G(H)] \oplus M^{c}(3,k),$ with dim $H_{0} = 12$,
\item $H_{0} = k[G(H)] \oplus M^{c}(2,k)^{3},$ with dim $H_{0} = 15$,
\item $H_{0} = k[G(H)] \oplus M^{c}(3,k)^{2},$ with dim $H_{0} = 21$,
\item $H_{0} = k[G(H)] \oplus M^{c}(2,k)^{3} \oplus M^{c}(3,k) ,$ 
with dim $H_{0} = 24$.
\end{enumerate}

Since all skew primitive elements of $H$ are trivial, by \cite[Cor. 4.3]{BD} 
we have that 
\begin{equation}\label{cuentas}
27 = \dim H > \dim H_{1} \geq (1 + 2n_{1})3 + \sum_{i=1}^{t}n_{i}^{2}.
\end{equation}

\noindent Replacing the corresponding dimensions in equation (\ref{cuentas}), 
it follows that no one of 
the cases (1), \ldots , (4) is possible, obtaining in this way a 
contradiction to our assumption.
\end{proof}

\begin{center}
{\sc Acknowledgment}

\end{center}
I thank Prof. N. Andruskiewitsch and Prof. 
H.-J. Schneider for their excellent guidance and for fruitful discussions
which helped to simplify and clarify the exposition. My thanks also to 
the referee for his or her suggestions to improve the proof of Theorem 
\ref{dim1}.
 
\bibliographystyle{amsbeta}

\end{document}